\documentclass[12pt]{article}
\usepackage{latexsym,amsfonts,amssymb,amsmath}
\usepackage{color}

\usepackage{colordvi,multicol}

\newtheorem{thm}{Theorem}[section]
\newtheorem{cor}[thm]{Corollary}

\newtheorem{rem}[thm]{Remark}

\date{}
\begin{document}

\title{\bf Quantitative Versions of the Two-dimensional  Gaussian Product Inequalities  }
 \author{ Ze-Chun Hu, Han Zhao, Qian-Qian Zhou\footnote{Corresponding author}\\ \\
 {\small College of Mathematics, Sichuan  University,
 Chengdu 610065, China}\\
 {\small zchu@scu.edu.cn}\\ \\
 {\small College of Mathematics, Sichuan  University,
 Chengdu 610065, China}\\
{\small 1339875802@qq.com}\\ \\
 {\small School of Science, Nanjing University of Posts and Telecommunications, Nanjing  210023, China}\\
 {\small qianqzhou@yeah.net}}
\maketitle

\begin{abstract}

\noindent  The Gaussian product inequality (GPI) conjecture is one of the most  famous  inequalities associated with Gaussian distributions and has attracted a lot of concerns. In this note, we investigate the quantitative versions of the two-dimensional Gaussian product inequalities. For any centered non-degenerate two-dimensional Gaussian random vector $(X_1, X_2)$  with variances $\sigma_1^2, \sigma_2^2$ and the correlation coefficient $\rho$,   we  prove that for any real numbers $\alpha_1, \alpha_2\in (-1,0)$ or $\alpha_1, \alpha_2\in (0,\infty)$,  it holds that
$${\bf E}[|X_1|^{\alpha_1}|X_2|^{\alpha_2}]-{\bf E}[|X_1|^{\alpha_1}]{\bf E}[|X_2|^{\alpha_2}]\ge f(\sigma_1,\sigma_2,\alpha_1, \alpha_2, \rho)\ge 0,
$$
where the function $f(\sigma_1,\sigma_2,\alpha_1, \alpha_2, \rho)$ will be given explicitly by Gamma function and is positive when $\rho\neq 0$. When $-1<\alpha_1<0$ and $\alpha_2>0,$ Russell and Sun (arXiv: 2205.10231v1) proved the ``opposite Gaussian product inequality", of which we will also give a quantitative version. These quantitative  inequalities are derived by employing the hypergeometric functions and the generalized hypergeometric functions.

\end{abstract}

\noindent  {\it MSC:} Primary 60E15; secondary 62H12

\noindent  {\it Keywords:} Gaussian product inequality conjecture, Quantitative inequality, Hypergeometric function.

%

\section{Introduction}



The  Gaussian product inequality (GPI) conjecture  is one of the most  famous  inequalities associated with Gaussian distributions, which states  that for any $d$-dimensional real-valued centered Gaussian random vector $(X_1, \ldots, X_d)$,
\begin{gather}\label{a1}
{\bf E}\left[\prod_{j=1}^{d}X_j^{2m}\right]\ge \prod_{j=1}^d{\bf E}\left[X_j^{2m}\right],\ m\in \mathbb{N}.
\end{gather}
In \cite{LW12}, Li and Wei proposed the following improved version of the GPI conjecture:
\begin{gather}\label{a2}
{\bf E}\left[\prod_{j=1}^{d}|X_j|^{\alpha_j}\right]\ge \prod_{j=1}^d{\bf E}\left[|X_j|^{\alpha_j}\right],
\end{gather}
where $\alpha_j,j=1, 2,\ldots, d$, are nonnegative real numbers.

Up to now, the GPI  \eqref{a1}  and the GPI \eqref{a2} are still open, however, some special cases have been proved by using various tools.  Frenkel \cite{F08} proved \eqref{a1} with $m=1$ (or \eqref{a2} for the case $\alpha_j=2$)  by using  algebraic methods.  Malicet et al. \cite{MNPP16} gave an analytic proof to Frenkel's result among other things.  Wei \cite{W14} proved a stronger version of \eqref{a2} for $\alpha_j\in (-1,0)$  as follows:
\begin{eqnarray}\label{a3}
{\bf E}\left[\prod_{j=1}^{d}|X_j|^{\alpha_j}\right]\ge {\bf E}\left[\prod_{j=1}^k|X_j|^{\alpha_j}\right]{\bf E}\left[\prod_{j=k+1}^d|X_j|^{\alpha_j}\right],\ \forall 1\leq k\leq d-1.
\end{eqnarray}

By Karlin and Rinott \cite{KR81}, we know that \eqref{a1}  and \eqref{a2} hold for $\mathbf{X}=(X_1,\ldots,X_d)$ when the density of  $|{\bf X}|=(|X_1|,|X_2|,\ldots, |X_d|)$ satisfies the so-called condition {\bf MTP$_2$}, and for any non-degenerate 2-dimensional centered Gaussian random vector $(X_1, X_2),$ $(|X_1|,|X_2|)$  has an {\bf MTP$_2$} density. Thus \eqref{a1}  and \eqref{a2} hold for $d=2.$

Lan et al. \cite{LHW20} used the hypergeometric functions to prove the following 3-dimensional GPI: for any $m_1, m_2\in \mathbb{N}$ and any centered Gaussian random vector $(X_1, X_2, X_3)$,
\begin{eqnarray}\label{a4}
{\bf E}\left[X_1^{2m_1}X_2^{2m_2}X_3^{2m_2}\right]\ge {\bf E}\left[X_1^{2m_1}\right]{\bf E}\left[X_2^{2m_2}\right]{\bf E}\left[X_3^{2m_2}\right].
\end{eqnarray}

Genest and Ouimet \cite{GO22} proved that if there exists a matrix $C\in [0, +\infty)^{d\times d}$ such that\linebreak  $(X_1, X_2, \ldots, X_d)=(Z_1, Z_2, \ldots, Z_d)C$ in law, where $(Z_1, Z_2, \ldots, Z_d)$ is a $d$-dimensional standard Gaussian random vector, the following stronger version of \eqref{a1} holds:
\begin{eqnarray}\label{a1-1}
{\bf E}\left[\prod_{j=1}^{d}X_j^{2m_j}\right]\ge \mathbf{E}\left[\prod_{j=1}^kX_j^{2m_j}\right]\mathbf{E}
\left[\prod_{j=k+1}^nX_j^{2m_j}\right],\ m_j\in \mathbb{N},j=1,\ldots,n,\ \forall 1\leq k\leq n-1.
\end{eqnarray}
Russell and Sun \cite{RS22-a} proved among other things that \eqref{a1-1} holds if all the correlation coefficients are nonnegative. Edelmann et al. \cite{ERR22} extended \eqref{a1-1} to the multivarite gamma distributions. As to other related work, we refer to Russell and Sun \cite{RS22-b}, Russell and Sun \cite{RS22-c}, Genest and Ouiment \cite {GO22-b}.

Recently, De et al. \cite{NNS20} obtained a range of  quantitative correlation inequalities, including the  quantitative version of the Gaussian correlation inequality   proved by Royen \cite{R14}, and the quantitative version of the well-known  Fortuin-Kasteleyn-Ginibre  (FKG) inequality for monotone functions over
any finite product probability space. Motivated by  \cite{NNS20}, we will  consider the quantitative versions of the two-dimensional Gaussian product inequalities in this note.

For any centered non-degenerate two-dimensional Gaussian random vector $(X_1, X_2)$  with variances $\sigma_1^2, \sigma_2^2$ and the correlation coefficient $\rho$,   we will prove that for any real numbers $\alpha_1, \alpha_2\in (-1,0)$ or $\alpha_1, \alpha_2\in (0,\infty)$,  it holds that
$${\bf E}[|X_1|^{\alpha_1}|X_2|^{\alpha_2}]-{\bf E}[|X_1|^{\alpha_1}]{\bf E}[|X_2|^{\alpha_2}]\ge f(\sigma_1,\sigma_2,\alpha_1, \alpha_2, \rho)\ge 0,
$$
where the function $f(\sigma_1,\sigma_2,\alpha_1, \alpha_2, \rho)$ will be given explicitly by Gamma function and is positive when $\rho\neq 0$. When $-1<\alpha_1<0$ and $\alpha_2>0,$ Russell and Sun \cite{RS22-c} proved the ``opposite Gaussian product inequality", of which we will also give a quantitative version. These quantitative  inequalities are derived by employing the hypergeometric functions and the generalized hypergeometric functions..

The rest of this note is organized as follows. In Section 2,  we  present the main results and some corollaries. In Section 3, we give  the proofs. In the final section, we give some remarks.

\section{Main results}\setcounter{equation}{0}
Throughout this note, any Gaussian random variable is assumed to be real-valued and non-degenerate. Our main results are as follows.

\begin{thm}\label{th-1}
Let $(X_1,X_2)$ be   centered bivariate  Gaussian random variables  with  $\textbf{E}[X_1^{2}]={\sigma_1}^2$, $\textbf{E}[X_2^{2}]={\sigma_2}^2$ and the correlation coefficient $\rho$ satifying $|\rho|<1$. Then for any real numbers $\alpha_1, \alpha_2\in (-1,0)$ or $\alpha_1, \alpha_2\in (0,\infty)$,
\begin{align}\label{th-1-a}
	\textbf{E}[{ \vert X_1 \vert }^{\alpha_1}{ \vert X_2 \vert }^{\alpha_2}]-\textbf{E}[{ \vert X_1 \vert }^{\alpha_1}]\textbf{E}[{ \vert X_2 \vert }^{\alpha_2}]\geq f(\sigma_1,\sigma_2,\alpha_1, \alpha_2, \rho),
\end{align}
where the nonnegative function $f(\cdot)$ is defined by
\begin{eqnarray*}
&&f(\sigma_1,\sigma_2,\alpha_1, \alpha_2, \rho)\\
&&:=\left\{
		\begin{array}{ll}			\frac{2^{\frac{\alpha_1+\alpha_2}{2}}\alpha_1\alpha_2\sigma_1^{\alpha_1}\sigma_2^{\alpha_2}\rho^2}{2\pi}
\Gamma(\frac{\alpha_1+1}{2})\Gamma(\frac{\alpha_2+1}{2}),&  \mbox{if}\ -1<\alpha_1,\alpha_2<0\ \mbox{or}\ 0<\alpha_1,\alpha_2\leq 2\ \mbox{or}\  \alpha_1, \alpha_2>2,\\
\frac{2^{\frac{\alpha_1+\alpha_2}{2}}\alpha_1\alpha_2\sigma_1^{\alpha_1}
\sigma_2^{\alpha_2}\rho^2}{4\sqrt{\pi}}\Gamma(\frac{\alpha_1+\alpha_2-1}{2}),
&\mbox{if}\ \alpha_1>2,0<\alpha_2<2\ \mbox{or}\ 0<\alpha_1<2, \alpha_2>2.
		\end{array} \right.
	\end{eqnarray*}
\end{thm}

When the real numbers $\alpha_1 $  and $\alpha_2$ in Theorem \ref{th-1} have opposite signs, Russell and Sun \cite{RS22-c} proved the following opposite GPI:
\begin{eqnarray}\label{2.2}
	\textbf{E}[{ \vert X_1 \vert }^{\alpha_1}{ \vert X_2 \vert }^{\alpha_2}]\leq \textbf{E}[{ \vert X_1 \vert }^{\alpha_1}]\textbf{E}[{ \vert X_2 \vert }^{\alpha_2}].
\end{eqnarray}
  We have the  quantitative version of the GPI (\ref{2.2}) as follows.

\begin{thm}\label{thm-2} Let $(X_1,X_2)$ be   centered bivariate  Gaussian random variables  with  $\textbf{E}[X_1^{2}]={\sigma_1}^2$, $\textbf{E}[X_2^{2}]={\sigma_2}^2$ and the correlation coefficient $\rho$ satifying $|\rho|<1$. Let $\alpha_1\in (-1,0)$ and $\alpha_2\in (0,\infty)$.

(i) When $0<\alpha_2\le2$,
	\begin{eqnarray}\label{thm-2-a}		&&\frac{2^{\frac{\alpha_1+\alpha_2}{2}}\alpha_1\alpha_2\sigma_1^{\alpha_1}
\sigma_2^{\alpha_2}\rho^2}{2\pi}\Gamma\left(\frac{\alpha_1+1}{2}\right)
\Gamma\left(\frac{\alpha_2+1}{2}\right)F\left(1-\frac{\alpha_1}{2},1-\frac{\alpha_2}{2};
\frac{3}{2};1\right)\nonumber\\
&&\leq \textbf{E}\left[{ \vert X_1 \vert }^{\alpha_1}{ \vert X_2 \vert }^{\alpha_2}\right]-\textbf{E}\left[{ \vert X_1 \vert }^{\alpha_1}\right]\textbf{E}\left[{ \vert X_2 \vert}^{\alpha_2}\right]\nonumber\\
&&\leq		\frac{2^{\frac{\alpha_1+\alpha_2}{2}}\alpha_1\alpha_2\sigma_1^{\alpha_1}
\sigma_2^{\alpha_2}\rho^2}{2\pi}\Gamma\left(\frac{\alpha_1+1}{2}\right)
\Gamma\left(\frac{\alpha_2+1}{2}\right)
		\leq0,
	\end{eqnarray}
where $F(\cdot)$  is the  hypergeometric function defined by
\begin{eqnarray}\label{GHF}
	F(a,b;c;z):=\sum\limits_{n=0}^{+\infty}\frac{(a)_n(b)_n}{(c)_n}\cdot\frac{z^n}{n!},\  |z|\leq 1,
\end{eqnarray}
and, for $\alpha\in\mathbb{R}$,
\begin{align*}
	(\alpha)_n:=\left\{
	\begin{array}{cl}
		\alpha(\alpha+1)\dots(\alpha+n-1), & n  \geq  1,\\
		1, & n=0,\alpha\neq 0.\\
	\end{array} \right.
\end{align*}

(ii) When $\alpha_2>2$,
	\begin{eqnarray}\label{thm-2-b}			&&\frac{2^{\frac{\alpha_1+\alpha_2}{2}}\alpha_1\alpha_2\sigma_1^{\alpha_1}\sigma_2^{\alpha_2}\rho^2}{2\pi}
\Gamma\left(\frac{\alpha_1+1}{2}\right)
\Gamma\left(\frac{\alpha_2+1}{2}\right)\nonumber\\
&&\leq \textbf{E}\left[{ \vert X_1 \vert }^{\alpha_1}{ \vert X_2 \vert }^{\alpha_2}\right]-\textbf{E}\left[{ \vert X_1 \vert }^{\alpha_1}\right]\textbf{E}\left[{ \vert X_2 \vert}^{\alpha_2}\right]\nonumber\\
&&\leq
		\min\left\{ \frac{2^{\frac{\alpha_1+\alpha_2}{2}}\alpha_1\alpha_2\sigma_1^{\alpha_1}\sigma_2^{\alpha_2}
\rho^2}{2\pi}\Gamma\left(\frac{\alpha_1+1}{2}\right)\Gamma\left(\frac{\alpha_2+1}{2}\right)
F\left(1-\frac{\alpha_1}{2},1-\frac{\alpha_2}{2};\frac{3}{2};1\right),\, 0\right\}.\quad\quad
	\end{eqnarray}
\end{thm}

By Theorem \ref{th-1}, we have the following corollaries.

\begin{cor}
	 Let $\alpha_1,\alpha_2\in \{1,2\}$  in   Theorem \ref{th-1}, then we have
\begin{eqnarray*}
&&\textbf{E}[{ \vert X_1 \vert }{ \vert X_2 \vert }]-\textbf{E}[{ \vert X_1 \vert }]\textbf{E}[{ \vert X_2 \vert }]\geq\frac{\sigma_1\sigma_2\rho^2}{\pi},\\
&&\textbf{E}[{ \vert X_1 \vert }{ X_2}^{2}]-\textbf{E}[{ \vert X_1 \vert }]\textbf{E}[{X_2}^{2}]\geq\frac{\sqrt{2}\sigma_1\sigma_2^{2}\rho^2}{\sqrt{\pi}},\\
&&\textbf{E}[{ X_1 }^{2}{X_2}^{2}]-\textbf{E}[{X_1}^{2}]\textbf{E}[{ X_2 }^{2}]\geq 2\sigma_1^2\sigma^2_2\rho^2.
\end{eqnarray*}
\end{cor}

\begin{cor}
Under the conditions of  Theorem \ref{th-1}, suppose that  $\alpha_2=1$, $\alpha_1$ is an  integer satisfying  $\alpha_1=m>2$. Then
	\begin{eqnarray*}
		\textbf{E}[{ \vert X_1 \vert }^{m}{ \vert X_2 \vert }]-\textbf{E}[{ \vert X_1 \vert }^{m}]\textbf{E}[{ \vert X_2 \vert }]\geq\left\{
		\begin{array}{rcl}
			\frac{(m-2)!!\,m\sigma_1^{m}\sigma_2\rho^{2}}{\sqrt{2\pi}}, & \text{if }  m\; is\; even,\\
			\frac{(m-2)!!\,m\sigma_1^{m}\sigma_2\rho^{2}}{2},&  \text{if }    m\;is\;odd.
		\end{array} \right.
	\end{eqnarray*}
\end{cor}

\begin{cor}
	Under the conditions of  Theorem \ref{th-1}, suppose that  $\alpha_1$ and $\alpha_2$ are integers and  $\alpha_1=m,\alpha_2=n$ with $m,n\in (2,\infty)$,

	(\expandafter{\romannumeral+1})  if  $m, n$ are both even integers, then
\begin{align*}
	\textbf{E}[{ X_1 }^{m}{X_2}^{n}]-\textbf{E}[{X_1}^{m}]\textbf{E}[{ X_2 }^{n}]&\geq\frac{(m-1)!!(n-1)!!mn\sigma_1^m\sigma_2^n\rho^2}{2};
\end{align*}

(\expandafter{\romannumeral+2}) if $m, n$ are both odd integers, then
\begin{align*}
	\textbf{E}[{ \vert X_1 \vert }^{m}{ \vert X_2 \vert }^{n}]-\textbf{E}[{ \vert X_1 \vert }^{m}]\textbf{E}[{ \vert X_2 \vert }^{n}]&\geq\frac{(m-1)!!(n-1)!!mn\sigma_1^m\sigma_2^n\rho^2}{\pi};
\end{align*}

(\expandafter{\romannumeral+3}) if one of $m$  and $n$   is odd and the other is even, then
\begin{align*}
	\textbf{E}[{ \vert X_1 \vert }^{m}{ \vert X_2 \vert }^{n}]-\textbf{E}[{ \vert X_1 \vert }^{m}]\textbf{E}[{ \vert X_2 \vert }^{n}]&\geq\frac{(m-1)!!(n-1)!!mn\sigma_1^m\sigma_2^n\rho^2}{\sqrt{2\pi}}.
\end{align*}
\end{cor}


\section{Proofs}


In this section, we give the proofs of Theorems \ref{th-1} and \ref{thm-2}.

\subsection{Proof of Theorem \ref{th-1}}
Obviously, we can assume that $\rho\neq 0$.  By the density function of centered Gaussian  random variable and the definition of Gamma function, we have that for $i=1,2$,
	\begin{eqnarray}\label{c0}
		\textbf{E}[{ \vert X_i \vert }^{\alpha_i}]&=&\int_{-\infty}^{+\infty}\frac{1}{\sqrt{2\pi}\sigma_i}{\vert x \vert}^{\alpha_i}e^{-\frac{x^2}{2\sigma_i^2}}dx=\frac{\sqrt{2}}{\sqrt{\pi}\sigma_i}
\int_{0}^{+\infty}x^{\alpha_i}e^{-\frac{x^2}{2\sigma_i^2}}dx\nonumber\\
		&=&\frac{2^{\frac{\alpha_i}{2}}\sigma_i^{\alpha_i}}{\sqrt{\pi}}
\int_{0}^{+\infty}y^{\frac{\alpha_i-1}{2}}e^{-y}dy=
\frac{2^{\frac{\alpha_i}{2}}\sigma_i^{\alpha_i}}{\sqrt{\pi}}
\Gamma\left(\frac{\alpha_i+1}{2}\right).
	\end{eqnarray}
Then
\begin{gather}\label{c1}
	\textbf{E}[{ \vert X_1 \vert }^{\alpha_1}]\textbf{E}[{ \vert X_2 \vert }^{\alpha_2}]=\frac{2^{\frac{\alpha_1+\alpha_2}{2}}
\sigma_1^{\alpha_1}\sigma_2^{\alpha_2}}{\pi}
\Gamma\left(\frac{\alpha_1+1}{2}\right)\Gamma\left(\frac{\alpha_2+1}{2}\right).
\end{gather}

Since $|\rho|<1,$ by Nabeya \cite{NS51}, we know that
\begin{gather}\label{c2}
	\textbf{E}[{ \vert X_1 \vert }^{\alpha_1}{ \vert X_2 \vert }^{\alpha_2}]=\frac{2^{\frac{\alpha_1+\alpha_2}{2}}
\sigma_1^{\alpha_1}\sigma_2^{\alpha_2}}{\pi}
\Gamma\left(\frac{\alpha_1+1}{2}\right)\Gamma\left(\frac{\alpha_2+1}{2}\right)
F\left(-\frac{\alpha_1}{2},-\frac{\alpha_2}{2};\frac{1}{2};\rho^2\right),
\end{gather}
where $F(\cdot)$  is the  hypergeometric function defined by (\ref{GHF}).

It follows from \eqref{c1}  and \eqref{c2} that
\begin{eqnarray}\label{c3}
&&\textbf{E}[{ \vert X_1 \vert }^{\alpha_1}{ \vert X_2 \vert }^{\alpha_2}]-\textbf{E}[{ \vert X_1 \vert }^{\alpha_1}]\textbf{E}[{ \vert X_2 \vert }^{\alpha_2}]\nonumber\\
&&=\frac{2^{\frac{\alpha_1+\alpha_2}{2}}\sigma_1^{\alpha_1}\sigma_2^{\alpha_2}}{\pi}
\Gamma\left(\frac{\alpha_1+1}{2}\right)\Gamma\left(\frac{\alpha_2+1}{2}\right)
\left[F\left(-\frac{\alpha_1}{2},-\frac{\alpha_2}{2};\frac{1}{2};\rho^2\right)-1\right].
\end{eqnarray}

Since $\frac{2^{\frac{\alpha_1+\alpha_2}{2}}\sigma_1^{\alpha_1}\sigma_2^{\alpha_2}}{\pi}
\Gamma(\frac{\alpha_1+1}{2})\Gamma(\frac{\alpha_2+1}{2})>0$ for any $-1<\alpha_1, \alpha_2<0$ or $\alpha_1, \alpha_2>0,$
it is enough to find the lower bound of
 $F(-\frac{\alpha_1}{2},-\frac{\alpha_2}{2};\frac{1}{2};\rho^2)-1.$
At first, we have that
\begin{align}\label{c4}
	F\left(-\frac{\alpha_1}{2},-\frac{\alpha_2}{2};\frac{1}{2};\rho^2\right)-1
&=\sum\limits_{k=1}^{+\infty}\frac{(-\frac{\alpha_1}{2})_k(-\frac{\alpha_2}{2})_k}
{(\frac{1}{2})_k}\cdot\frac{({\rho^2})^k}{k!}\nonumber\\
	&=\sum\limits_{k=0}^{+\infty}\frac{(-\frac{\alpha_1}{2})_{k+1}
(-\frac{\alpha_2}{2})_{k+1}}{(\frac{1}{2})_{k+1}}\cdot\frac{({\rho^2})^{k+1}}{(k+1)!}\nonumber\\
&=\frac{\rho^2\alpha_1\alpha_2}{2}\sum\limits_{k=0}^{+\infty}
\frac{(1-\frac{\alpha_1}{2})_{k}(1-\frac{\alpha_2}{2})_{k}(1)_{k}}
{(\frac{3}{2})_{k}(2)_{k}}\cdot\frac{{\rho^2}^k}{(k)!}\nonumber\\
	&=\frac{\rho^2\alpha_1\alpha_2}{2}\     {_{3}F_{2}}\left(1-\frac{\alpha_1}{2},1-\frac{\alpha_2}{2},1;\frac{3}{2},2;\rho^2\right),
\end{align}
where $_{p}F_{q}(a_1,\dots,a_p;b_1,\dots,b_q;z)$  is the generalized hypergeometric function defined by
\begin{align*}
	_{p}F_{q}(a_1,\dots,a_p;b_1,\dots,b_q;z):=\sum\limits_{k=0}^{+\infty}\frac{(a_1)_n\dots(a_p)_n}{(b_1)_n\dots(b_q)_n}\cdot\frac{z^n}{n!}.
\end{align*}
Especially,  when $p=2$ and $q=1$, then $_{p}F_{q}(a_1,\dots,a_p;b_1,\dots,b_q;z)=F(a_1,a_2;b_1;z)$.
By Andrews et al.  \cite[(2.2.2), P. 67]{AAR99},  we know that
\begin{eqnarray}\label{c5}
&&_{p+1}F_{q+1}(a_1,\dots,a_p,a_{p+1};b_1,\dots,b_q,b_{q+1};z)\nonumber\\
&&=\frac{\Gamma(b_{q+1})}{\Gamma(a_{p+1})\Gamma(b_{q+1}-a_{p+1})}
\int_{0}^{1}t^{a_{p+1}-1}(1-t)^{b_{q+1}-a_{p+1}-1}\,
_{p}F_{q}(a_1,\dots,a_p;b_1,\dots,b_q;zt)dt.\ \ \
\end{eqnarray}
By \eqref{c4}  and \eqref{c5}, we get
\begin{align}\label{c6}
	F\left(-\frac{\alpha_1}{2},-\frac{\alpha_2}{2};\frac{1}{2};\rho^2\right)-1
&=\frac{\rho^2\alpha_1\alpha_2}{2}\frac{\Gamma(2)}{\Gamma(1)\Gamma(1)}\int_{0}^{1}
F\left(1-\frac{\alpha_1}{2},1-\frac{\alpha_2}{2};\frac{3}{2};\rho^2t\right)dt\nonumber\\	&=\frac{\alpha_1\alpha_2}{2}\int_{0}^{\rho^2}
F\left(1-\frac{\alpha_1}{2},1-\frac{\alpha_2}{2};\frac{3}{2};h\right)dh.
\end{align}

Secondly, the following property of the hypergeometric function,
\begin{align*}
	\frac{d}{dz}F(a,b;c;z)=\frac{ab}{c}F(a+1,b+1;c+1;z),
\end{align*}
(see Andrews et al. \cite[(2.5.1), P. 94]{AAR99}) and the Euler transformation (see  Andrews et al. \cite[Theorem 2.2.5]{AAR99} or Rainville \cite[Chapter 4, Theorem 21]{R60}) imply that
 \begin{eqnarray}\label{c7}
&&\frac{d}{dh}F\left(1-\frac{\alpha_1}{2},1-\frac{\alpha_2}{2};\frac{3}{2};h\right)\nonumber\\
&&=\frac{2}{3}\left(1-\frac{\alpha_1}{2}\right)\left(1-\frac{\alpha_2}{2}\right)
F\left(2-\frac{\alpha_1}{2},2-\frac{\alpha_2}{2};\frac{5}{2};h\right)\nonumber\\	&&=\frac{2}{3}\left(1-\frac{\alpha_1}{2}\right)\left(1-\frac{\alpha_2}{2}\right)
(1-h)^{\frac{1+\alpha_1+\alpha_2}{2}-2}
F\left(\frac{1}{2}+\frac{\alpha_1}{2},\frac{1}{2}+\frac{\alpha_2}{2};\frac{5}{2};h\right).
\end{eqnarray}
Then, in order  to  estimate $F(-\frac{\alpha_1}{2},-\frac{\alpha_2}{2};\frac{1}{2};\rho^2)-1$, we  divide it  into the following two cases.

\textbf{Case I}:  When  $-1<\alpha_1,\alpha_2<0$ or $0<\alpha_1,\alpha_2\leq2,$ or $\alpha_1, \alpha_2>2.$

In this case,  $(1-\frac{\alpha_1}{2}) (1-\frac{\alpha_2}{2} )\geq0$ and $\alpha_1\alpha_2>0.$ Thus, the derivative of $\frac{d}{dh}F(1-\frac{\alpha_1}{2},1-\frac{\alpha_2}{2};\frac{3}{2};h)$ in \eqref{c7} is nonnegative. Therefore, by \eqref{c6} we obtain that
 \begin{align*}
	F\left(-\frac{\alpha_1}{2},-\frac{\alpha_2}{2};\frac{1}{2};\rho^2\right)-1&
\geq \frac{\alpha_1\alpha_2}{2}\int_{0}^{\rho^2}
F\left(1-\frac{\alpha_1}{2},1-\frac{\alpha_2}{2};\frac{3}{2};0\right)dh
=\frac{\alpha_1\alpha_2}{2}\cdot\rho^2.
\end{align*}
By (\ref{c3}), we obtain
\begin{align*}
	\textbf{E}[{ \vert X_1 \vert }^{\alpha_1}{ \vert X_2 \vert }^{\alpha_2}]-\textbf{E}[{ \vert X_1 \vert }^{\alpha_1}]\textbf{E}[{ \vert X_2 \vert }^{\alpha_2}]&\geq\frac{2^{\frac{\alpha_1+\alpha_2}{2}}
\sigma_1^{\alpha_1}\sigma_2^{\alpha_2}}{\pi}
\Gamma\left(\frac{\alpha_1+1}{2}\right)\Gamma\left(\frac{\alpha_2+1}{2}\right)
\frac{\alpha_1\alpha_2}{2}\cdot\rho^2\\
	&=\frac{2^{\frac{\alpha_1+\alpha_2}{2}}\alpha_1\alpha_2\sigma_1^{\alpha_1}
\sigma_2^{\alpha_2}\rho^2}{2\pi}\Gamma\left(\frac{\alpha_1+1}{2}\right)
\Gamma\left(\frac{\alpha_2+1}{2}\right).
\end{align*}

\textbf{Case II}:
 When $2<\alpha_1$ and $0<\alpha_2<2$ or $0<\alpha_1<2$ and $2<\alpha_2$.

Without loss of generality, we assume  that $2<\alpha_1$ and $0<\alpha_2<2$.   Then  $(1-\frac{\alpha_1}{2})(1-\frac{\alpha_2}{2})<0$  and $\alpha_1\alpha_2>0.$  Thus, in this case, by (\ref{c7}), we know that
\begin{align*}
	\frac{d}{dh}F\left(1-\frac{\alpha_1}{2},1-\frac{\alpha_2}{2};\frac{3}{2};h\right)<0.
\end{align*}
Then   $F(1-\frac{\alpha_1}{2},1-\frac{\alpha_2}{2};\frac{3}{2};h)$ reaches its minimum at $h=1$,
 and by Rainville \cite[Chapter 4, Theorem 18]{R60}, the minimum value is
\begin{align*}
	F\left(1-\frac{\alpha_1}{2},1-\frac{\alpha_2}{2};\frac{3}{2};1\right)
=\frac{\Gamma(\frac{3}{2})\Gamma(\frac{\alpha_1+\alpha_2}{2}-\frac{1}{2})}{\Gamma(\frac{\alpha_1}{2}+\frac{1}{2})\Gamma(\frac{\alpha_2}{2}+\frac{1}{2})}.
\end{align*}
Then by (\ref{c6}), we get
\begin{align*}
	F\left(-\frac{\alpha_1}{2},-\frac{\alpha_2}{2};\frac{1}{2};\rho^2\right)-1
	&\geq\frac{\alpha_1\alpha_2}{2}\cdot\rho^2\frac{\Gamma(\frac{3}{2})\Gamma(\frac{\alpha_1+\alpha_2}{2}-\frac{1}{2})}{\Gamma(\frac{\alpha_1}{2}+\frac{1}{2})\Gamma(\frac{\alpha_2}{2}+\frac{1}{2})}.
\end{align*}
By (\ref{c3}), we obtain
\begin{eqnarray*}
&&\textbf{E}[{ \vert X_1 \vert }^{\alpha_1}{ \vert X_2 \vert }^{\alpha_2}]-\textbf{E}[{ \vert X_1 \vert }^{\alpha_1}]\textbf{E}[{ \vert X_2 \vert }^{\alpha_2}]\\
&&\geq\frac{2^{\frac{\alpha_1+\alpha_2}{2}}\sigma_1^{\alpha_1}\sigma_2^{\alpha_2}}{\pi}
\Gamma\left(\frac{\alpha_1+1}{2}\right)
\Gamma\left(\frac{\alpha_2+1}{2}\right)\frac{\alpha_1\alpha_2}{2}\cdot\rho^2
\frac{\Gamma(\frac{3}{2})\Gamma(\frac{\alpha_1+\alpha_2}{2}-\frac{1}{2})}
{\Gamma(\frac{\alpha_1}{2}+\frac{1}{2})\Gamma(\frac{\alpha_2}{2}+\frac{1}{2})}\\	&&=\frac{2^{\frac{\alpha_1+\alpha_2}{2}}\alpha_1\alpha_2
\sigma_1^{\alpha_1}\sigma_2^{\alpha_2}\rho^2}{4\sqrt{\pi}}
\Gamma\left(\frac{\alpha_1+\alpha_2-1}{2}\right).
\end{eqnarray*}

The proof is complete. \hfill\fbox

\subsection{Proof of Theorem  \ref{thm-2}}

By the proof   of Theorem \ref{th-1}, we know that
\begin{eqnarray}
&&\textbf{E}[{ \vert X_1 \vert }^{\alpha_1}{ \vert X_2 \vert }^{\alpha_2}]-\textbf{E}[{ \vert X_1 \vert }^{\alpha_1}]\textbf{E}[{ \vert X_2 \vert }^{\alpha_2}]\nonumber\\
&&=\frac{2^{\frac{\alpha_1+\alpha_2}{2}}
\sigma_1^{\alpha_1}\sigma_2^{\alpha_2}}{\pi}
\Gamma\left(\frac{\alpha_1+1}{2}\right)
\Gamma\left(\frac{\alpha_2+1}{2}\right)
\left[F\left(-\frac{\alpha_1}{2},-\frac{\alpha_2}{2};\frac{1}{2};\rho^2\right)-1\right],
\label{Proof-thm-2-a}\\
&&F\left(-\frac{\alpha_1}{2},-\frac{\alpha_2}{2};\frac{1}{2};\rho^2\right)-1
=\frac{\alpha_1\alpha_2}{2}\int_{0}^{\rho^2}
F\left(1-\frac{\alpha_1}{2},1-\frac{\alpha_2}{2};\frac{3}{2};h\right)dh,\label{Proof-thm-2-b}
\end{eqnarray}
and
\begin{eqnarray*}	\frac{d}{dh}F\left(1-\frac{\alpha_1}{2},1-\frac{\alpha_2}{2};\frac{3}{2};h\right)
=\frac{2}{3}\left(1-\frac{\alpha_1}{2}\right)\left(1-\frac{\alpha_2}{2}\right)
(1-h)^{\frac{1+\alpha_1+\alpha_2}{2}-2}
F\left(\frac{1}{2}+\frac{\alpha_1}{2},\frac{1}{2}+\frac{\alpha_2}{2};\frac{5}{2};h\right).
\end{eqnarray*}

\textbf{Case I}: $0<\alpha_2\le 2$.  In this case $1-\frac{\alpha_2}{2}\geq 0$   and  thus
	\begin{align*}
		\frac{d}{dh}F\left(1-\frac{\alpha_1}{2},1-\frac{\alpha_2}{2};\frac{3}{2};h\right)&\geq0.
	\end{align*}
Then  by (\ref{Proof-thm-2-b}),  we get
	\begin{eqnarray*}
		\frac{\alpha_1\alpha_2\rho^2}{2}
F\left(1-\frac{\alpha_1}{2},1-\frac{\alpha_2}{2};\frac{3}{2};1\right)\leq F\left(-\frac{\alpha_1}{2},-\frac{\alpha_2}{2};\frac{1}{2};\rho^2\right)-1
\leq\frac{\alpha_1\alpha_2\rho^2}{2}\le 0,
	\end{eqnarray*}
which together with (\ref{Proof-thm-2-a}) implies (\ref{thm-2-a}).


	\textbf{Case II}: $\alpha_2>2$. In this case, $1-\frac{\alpha_2}{2}<0$, and thus
	\begin{align*}
		\frac{d}{dh}F\left(1-\frac{\alpha_1}{2},1-\frac{\alpha_2}{2};\frac{3}{2};h\right)&\leq0.
	\end{align*}
Then  by (\ref{Proof-thm-2-b}),  we get
	\begin{eqnarray*}
		\frac{\alpha_1\alpha_2\rho^2}{2}\leq F\left(-\frac{\alpha_1}{2},-\frac{\alpha_2}{2};\frac{1}{2};\rho^2\right)-1
\leq\frac{\alpha_1\alpha_2\rho^2}{2}
F\left(1-\frac{\alpha_1}{2},1-\frac{\alpha_2}{2};\frac{3}{2};1\right),
	\end{eqnarray*}
which together with (\ref{Proof-thm-2-a}) implies (\ref{thm-2-b}). 	

The proof is  complete.\hfill\fbox

\section{Remarks}

\begin{rem}\label{rem-4.1}
If $\rho$ in Theorem \ref{th-1} satisfies that $|\rho|=1$, then $X_2=\pm X_1$ in law and thus $\mathbf{E}[|X_1|^{\alpha_1}|X_2|^{\alpha_2}]=\mathbf{E}[|X_1|^{\alpha_1+\alpha_2}],
\sigma_1=\sigma_2$.
By (\ref{c0}) and $\Gamma(\frac{1}{2})=\sqrt{\pi}$, we know that
\begin{eqnarray*}
&&\mathbf{E}[|X_1|^{\alpha_1}|X_2|^{\alpha_2}]-\mathbf{E}[|X_1|^{\alpha_1}|]
\mathbf{E}[|X_2|^{\alpha_2}]\\
&&=\mathbf{E}[|X_1|^{\alpha_1+\alpha_2}]-\mathbf{E}[|X_1|^{\alpha_1}|]
\mathbf{E}[|X_2|^{\alpha_2}]\\
&&=\frac{2^{\frac{\alpha_1+\alpha_2}{2}}\sigma_1^{\alpha_1+\alpha_2}}{\sqrt{\pi}}
\Gamma\left(\frac{\alpha_1+\alpha_2+1}{2}\right)-\frac{2^{\frac{\alpha_1+\alpha_2}{2}}
\sigma_1^{\alpha_1+\alpha_2}}{\pi}
\Gamma\left(\frac{\alpha_1+1}{2}\right)\Gamma\left(\frac{\alpha_2+1}{2}\right)\\
&&=\frac{2^{\frac{\alpha_1+\alpha_2}{2}}\sigma_1^{\alpha_1}
\sigma_2^{\alpha_2}}{\pi}\Gamma\left(\frac{\alpha_1+1}{2}\right)
\Gamma\left(\frac{\alpha_2+1}{2}\right)\left[ \frac{\Gamma(\frac{\alpha_1+\alpha_2+1}{2})\Gamma(\frac{1}{2})}{\Gamma(\frac{\alpha_1+1}{2})
\Gamma(\frac{\alpha_2+1}{2})} -1 \right]\\
&&= \frac{2^{\frac{\alpha_1+\alpha_2}{2}}\sigma_1^{\alpha_1}\sigma_2^{\alpha_2}}{\pi}
\Gamma\left(\frac{\alpha_1+1}{2}\right)\Gamma\left(\frac{\alpha_2+1}{2}\right)\left[
F\left(-\frac{\alpha_1}{2}, -\frac{\alpha_2}{2};\frac{1}{2};1\right)-1\right],
\end{eqnarray*}
the last equality follows from  Rainville \cite[Chapter 4, Theorem 18]{R60}.
That is to say, (\ref{c3}) still holds in this case. Hence the result of Theorem \ref{th-1} still holds if $|\rho|=1$ and one of the following two conditions holds:

(i) $\alpha_1,\alpha_2\in (0,\infty)$;

(ii) $\alpha_1,\alpha_2\in (-1,0)$ and $\alpha_1+\alpha_2\in (-1,0)$.
\end{rem}

\begin{rem}\label{rem-4.2}
If $\rho$ in Theorem \ref{thm-2} satisfies that $|\rho|=1$, then by Remark \ref{rem-4.1} and the proof of Theorem \ref{thm-2} we know that the result of Theorem \ref{thm-2} still holds in this case.
\end{rem}

\begin{rem}\label{rem-4.3}
 A natural question is:

 How about the quantitative version of the 3-dimensional Gaussian product inequality (\ref{a4})?

We explored it, but we have obtained only some partial results up to now.
\end{rem}

\bigskip

{ \noindent {\bf\large Acknowledgments}  This work was supported
by the National Natural Science Foundation of China (12171335), the Science Development Project of Sichuan University (2020SCUNL201), and  the Scientific  Foundation of  Nanjing University of Posts and Telecommunications (NY221026).

\end{document}